\documentclass{amsart}

  \usepackage[all]{xy}

  \usepackage{epsf,epsfig,amsfonts,graphicx,color}

\usepackage{amsmath,amssymb}

\usepackage{eqparbox}

\usepackage{MnSymbol}

  \usepackage{url}
   
 \def\R{\mathbb R}

\def\C{\mathbb C}

\def\N {\mathbb N}

\def \nbhd { neighborhood }
\def \nbhds { neighborhoods }

\newcommand{\beq}{\begin{equation}}
\newcommand{\eeq}{\end{equation}}
\newcommand{\Leq}[1]{\label{#1}\end{equation}}

\newtheorem{theorem}{Theorem}[section]
\newtheorem{proposition}{Proposition}[section]
\newtheorem{corollary}{Corollary}[section]

\newtheorem{definition}{Definition}[section]
\newtheorem{remark}{Remark}[section]
\newtheorem{example}{Example}[section]

\newtheorem{working hyp}{Working Hypothesis}
\newtheorem{rmk}{Remark}

\def \Pc {Puiseux characteristic }

\begin{document}

\title{The Honest  Embedding Dimension of a Numerical Semigroup.}

\begin{abstract}  Attached    to a singular  analytic curve germ in $d$-space is  a  numerical
semigroup: a subset  $S$ of the non-negative integers  which is  closed under addition and  whose complement is
finite. 
 Conversely, associated to any numerical semigroup $S$ is a canonical mononial curve in $e$-space where $e$ is the number of minimal generators of the semigroup.    It may happen that  $d < e = e(S)$
 where $S$ is the semigroup of the curve in $d$-space.  Define the minimal (or `honest') embedding
of a  numerical semigroup to be the smallest $d$ such that $S$ is realized by a curve in
$d$-space.    Problem: characterize  the numerical semigroups having minimal embedding
dimension $d$.  The answer is known for the case $d=2$ of planar curves
and  reviewed  in an Appendix to this paper.
The case $d =3$ of   the problem is  open. Our main result is
a characterization of the  multiplicity $4$ numerical semigroups whose minimal embedding dimension is $3$.  See figure 1.    
The motivation for this work came from thinking about Legendrian curve singularities. 
\end{abstract}

\author{Richard Montgomery}
\address{Mathematics Department\\ University of California, Santa Cruz\\
Santa Cruz CA 95064}
\email{rmont@ucsc.edu}

\date{March 1,  2024}

\maketitle


A numerical semigroup is   a subset  $S \subset \N$ of the natural numbers  $\N$
 closed under addition and whose complement is finite.    
 Invariantly attached to any curve singularity in $d$-space is a numerical semigroup.
 See section \ref{semigp of a curve}.  For example, if $m < n$ are relatively prime integers then the planar curve
  $x = t^m, y = t^n + a t^{n+1} + b t^{n+2} + \ldots$ has the semigroup $<m, n>$ 
 attached to it
  regardless of the higher order terms $a t^{n+1} + b t^{n+1} + \ldots$ etc.
  Here $<m , n> =  \{ k m + \ell n:  k, \ell \in \N  \}$ denotes the
  semigroup  generated by $m$ and $n$, i.e. all sums of $m$ and $n$.    
   
Any numerical semigroup $S$ has a finite set 
$\{n_1, n_2, \ldots, n_k \} \subset S$ of generators,  meaning elements  such that every element of $S$ is a sum of the  $n_i$.
We write $S = <n_1, n_2, \ldots,  n_k >$.
Among all possible finite generating sets of $S$  the one with the 
smallest cardinality is unique.  The cardinality $e$ of this  minimal generating
set   $\{n_1, n_2, \ldots , n_e \}$  is called the embedding dimension of $S$. 
Underlying this `embedding' terminology   is a construction. 
Form the monomial curve $x_i = t^{n_i}, i=1, 2, \ldots, e$ in   $e$-space using the minimal generating set as exponents.
The semigroup attached to the monomial  curve is $S$.   But many    curves   besides the  monomial curve 
will have the same  semigroup $S$ attached to them.   Some   might even  lie  in a space of dimension  $d < e$.  
\begin{definition}
The minimal embedding dimension of a numerical semigroup $S$
is the minimal dimension $d$ such that the semigroup is that
of some analytic curve germ $(\C, 0) \to (\C^d,  0)$
\end{definition}

\begin{example} 
The semigroup $S = <4, 6, 13>$ consists of all integers which are sums of $4, 6$ and $13$,
these being its minimal generating set.  Its embedding dimension is $3$. Its canonical curve
is $x =t^4, y = t^6,  z = t^{13}$.
The plane curve  $x = t^4, y =t^6 + t^7$  has  this same
semigroup. The generator $13$ in the semigroup of this curve arises
 as the order of polynomial $p(x,y) = y^2 - x^3$,  
 when  pulled back to this plane curve :    $p(x,y) = 2 t^{13} + t^{14}$. 
So the minimal embedding dimension of $S$ is $2$. 
\end{example}

{\bf Problem.}  {\it Characterize  those numerical semigroups whose minimal embedding dimension
is } $d$. 

This problem has been solved for the case $d=2$ of plane curves.
See  Teissier \cite{Teissier}  proposition 3.2.1, on page 132. 
We recall and clarify  this  proposition  and sketch a proof based on the Puiseux characteristic
 in the   Appendix at the end of the paper.
The  case $d=3$ of space curves appears to be open.   See Castellanos \cite{Castellanos},
particularly   problem 2.4  and the examples in section 2 for  perspective.  

 Teissier states his  proposition 3.2.1  in a   convoluted way so as to hold for  all   $d \le e$. 
 In the appendix I   make what sense I can of his proposition as it applies to    $d > 2$.  
What he does is to  provide a list of   {\it sufficient} conditions amongst the $e$ minimal generators of a   semigroup  in order  for that semigroup to have minimal embedding dimension $d \le e$.  Teissier's  conditions  imply that  curves to which the  semigroup is attached are  complete intersections.  All plane curve branches are complete intersections
which allows his  conditions to be   necessary and sufficient when $d =2$.  
  Many space curve singularities fail to be complete intersections
and consequently Teissier's  conditions exclude  many  semigroups with minimial embedding dimension $3$.

The ``multiplicity'' $m$ 
of a numerical semigroup is its smallest nonzero element.
We have
$$me (S) \le e(S) \le m(S).$$
where $me(S), e(S)$, and $m(S)$
are the minimal emdedding dimension, embedding dimension and multiplicity
of the semigroup $S$.   (Refer to the first chapter of the book \cite{theBook}
for more standard terminology around numerical semigroups.)   If we  refine our problem   according to  
multiplicity  it becomes more   tractable. For example, semigroups of multiplicity
$2$ or $3$   have minimal embedding dimension less than or equal to $3$ by the above inequality
and so can be excluded as `trivial'  in the search for solutions to  the  problem for $d= 3$. 
In this paper  we will solve:  

\begin{figure}[h]
\scalebox{0.3}{\includegraphics{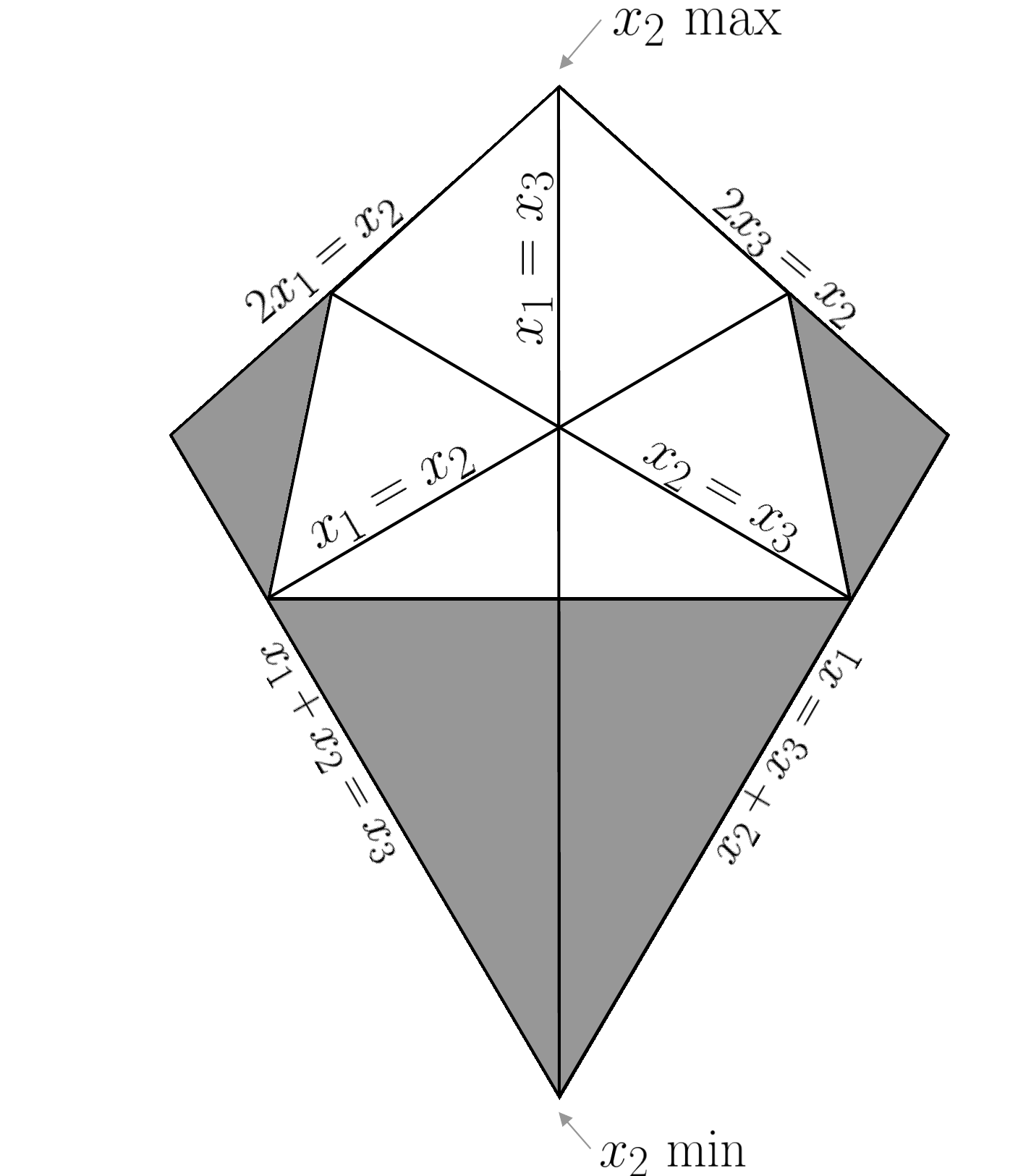}}
\caption{The Kunz cone for multiplicity $4$ is the  cone in $\R^3$ over the interior of this kite 
when placed on the plane $x_1 + x_2 + x_3 =1$. See equations (\ref{Kunzcone}).  
The set of  lattice points in and on the boundary of the cone for which $x_i \equiv i$ (mod 4) are in
bijection with multiplicity 4 numerical semigroups.  The set of such points lying in the   interior of the cone sweep out  the
semigroups with embedding dimension $4$.  Those having minimal embedding dimension $3$
correspond to the shaded region of the kite and the interiors of its edges.  Those having  minimal embedding $2$ correspond
to the left and right vertices of the kite.  {\bf Thanks to Emily O'Sullivan for the figure.}}  
\label{kunzkite}
\end{figure}

{\bf Problem.}  {\it Describe all   numerical semigroups 
whose  minimal embedding dimension  is $3$ and whose  multiplicity $4$.}

We solve this problem completely below.  See theorem \ref{thm: main} and  
figure \ref{kunzkite}.
To give the readers a taste of the solution,
represent an $m=4, e= 4$ numerical semigroup in the form 
$S = <4, n_1, n_2, n_3>$ with $4 < n_1 < n_2 < n_3$.  Since $e=4$ the four integers  $4, n_1, n_2, n_3$ must
represent all 4 congruence classes mod $4$. They also must satisfy
the  strict Kunz inequalities listed in  equation (\ref{Kunzcone}) below.   We will  show that  if
 $n_3$ is congruent to  $2$ mod $4$ then $S$ is not on our list.  So  for such $S$'s  we have $me(S) =4$.
 This fact corresponds to the white top triangle of the Kuntz kite of figure \ref{kunzkite}. 
And will see that  among those semigroups  of  the form $S = <4, 6, n_2, n_3>$   
if $me(S) = 3$ then  we must have  $n_3 = n_2 + 2$.

The Kunz cone associated to $m =4$ gives us a good language within which to solve the problem.  
We describe this  cone after recalling how to attach a semigroup to a curve.

\subsection{Motivation: Legendrian semigroups}

A Legendrian curve is a  space curve   tangent to a contact
 distribution in 3-space.  Legendrian curves are born from  plane curve
 singularities by a  process  known variously  as ``prolonging'',
 ``Nash blow-up'' or ``forming the conormal variety'' and which shares many properties with the classical blow-up
 of a singularity. 
 The  semigroup of a Legendrian curve has a  close  but poorly  understood relationship
with the semigroup of the planar curve which gave birth to it.     
    See for example \cite{Neto1},  \cite{ourBook}, or  \cite{zhit}.   The  following  questions   motivate  this work.  
     $$ \text{ What is  the set of   semigroups which arise from Legendrian curves? }$$
     \noindent More importantly, 
is this really an interesting question? For example, does the  semigroup of a singular Legendrian curve  encode interesting  unappreciated ``contact topological'' properties of the curve?  If the answer to the last question ends up being ``yes'' then the answer to the second question
would  also  be`yes'.

As a beginning  step towards answering these questions  we  have classified
the Legendrian semigroups of multiplicity 2, 3 and 4.  We hope to present our findings   in a companion paper.  
The first two cases (m = 2 and 3) are  easy and well-known.   Work on the  last case ($m =4$) naturally led into
the   problem solved here, since a Legendrian semigroup  must have minimal
embedding dimension 3 or less.

\section{The semigroup of a curve}
\label{semigp of a curve}

By a {\it curve}  we    mean a non-constant analytic map $c: (\C, 0) \to (\C^d, 0)$.
The notation $c: (\C, 0) \to (\C^d, 0)$ means that $c(0) = 0$.
It also indicates that our  true interest is  the germ of the curve,  meaning
its restriction to  any   small \nbhd of $0$.    
In order to attach  a numerical semigroup  to a curve   
we assume that the curve $c$  is well-parameterized:

\begin{definition}
A   curve is well-parameterized if it is one-to-one as a map
when restricted to a sufficiently small \nbhd of zero.  
\end{definition}
\noindent
If $c$ is not well-parameterized then it can be replaced
by another curve which is `the same' as $c$ and which is
well-parameterized. 
See the final subsection of this section. 

Recall the order of a power series converging in a \nbhd of  $t =0$.
The order of   $f(t) = \Sigma_{i \ge 0 } a_i t^i$, written ``$ord(f)$'', 
 is  the smallest $i$ such that $a_i \ne 0$.   For example $ord(t^3 + 7t^5) = 3$.
If  $\C\{t\}$ denotes the space of   power series which converge in some
\nbhd of $0$ then 
$$ord: \C\{t\} \to \N.$$
The order is a valuation, i.e. a semigroup homomorphism:  $ord(fg) = ord(f)+ ord(g)$.
Let $\mathcal O_0$ denote the ring of germs of
analytic functions $f: (\C^d, 0) \to \C$ defined in a \nbhd of $0$.
Let $\C \{x_1, \ldots , x_d \}$
denote the ring of convergent power series in the coordinates $x_i$ of $\C^d$,
converging in some \nbhd (depending on the series) of $0 \in \C^d$. 
 Pull-back defines a ring homomorphism  
$$c^* :  \C \{x_1, \ldots , x_d \} \to \C \{t\} 
$$
by sending $p \in \mathfrak m$ to $c^* p := p \circ c$.
Thus $c^*  \C \{x_1, \ldots , x_d \} \subset \C \{t \}$ is a subring.  The   semigroup of $c$
is the collection of all integers of the form $n = ord (c^* p)$
which   arise in this way.    In symbols  
\begin{definition}  [Semigroup of a curve.]  The semigroup
of the   analytic well-parameterized curve $c$ is the set of integers
$S = ord(c^*   \C \{x_1, \ldots , x_d \})$.  
\end{definition}
\noindent To see that $S$ is  closed under addition   use 
$ord(fg ) = ord(f) + ord(g)$ and $c^* (p q) = (c^*p) (c^* q)$.

\vskip .4cm

It may be helpful to write the pullback operation and order map  out in coordinates.  Write out both our curve 
and our function in coordinates: 
$$c(t) = (x_0 (t), x_1 (t), \ldots, x_d (t))$$ and $p = p(x_1, \ldots , x_d)$
where the $x_i$ are coordinates for $\C^d$.
Then $$f(t) = c^* p (t) = p(x_0 (t), x_1 (t), \ldots, x_d (t))$$
is an   analytic function of $t$.

\begin{rmk}  We could replace  ${\mathcal O}_0 = \C \{x_1, \ldots , x_d \}$
  by the ring  of polynomials   in the $x_i$
or the algebra of formal power series in the $x_i$
or even germs of smooth functions at $0$ and we would obtain the   same semigroup $S$
as $ord (c^* R)$.    This is because $c$ is well-parameterized and
as a consequence the   semigroup
$S$ is numerical and as such  has a ``conductor'' , a number $k$ such that all integers greater than or
equal to $k$ lie in $S$.
Considerations of degrees then show that by restricting oneself to polynomials in the $x_i$ of a fixed 
degree (roughly degree  $(k/m) +1 $ where $m$ is the multiplicity of the curve $c$) will suffice
to realize all elements of $S$.  
\end{rmk}

To see that our semigroup $S$   is ``numerical'', i.e. that  the  complement  of $S$ is finite,  
we must use  that  $S$ is well-parameterized.   
Suppose, by way of contradiction,  that   the complement
of $S$ were infinite.  Then the   g.c.d of
$S$ would not be $1$, but instead some integer $k >1$.
From this it follows that, roughly speaking,  all the exponents arising in all   power series
$c^* p,  p \in \mathfrak m$ would be divisible by $k$ and from this
that we could express   $c(t) = \gamma(t^k)$ for
some analytic curve $\gamma$.  This implies that   $c$
is not well-parameterized.    (We say `roughly speaking'
because we might need to reparameterize first:
$c(t) =  \gamma( \tau(t)^k)$ where $\tau'(0) \ne 0$.)

The reader may wish to refer to section 5 of Arnol'd \cite{Arnold} for 
a beautiful   perspective on the relation between a curve and its semigroup.   

\subsubsection{Analytic equivalence} Call two curves $c_1, c_2 : (\C, 0) \to (\C^d, 0)$
``analytically equivalent'' if there are germs of analytic
diffeomorphisms $\Psi: (\C^d, 0) \to (\C^d, 0)$ and
$\psi: (\C, 0) \to (\C, 0)$ such that
$c_2 = \Psi \circ c_1 \circ \phi$.    In the singularity literature
this equivalence relation is also called ``RL equivalence''.
We write $c_1 \sim c_2$ to mean that the two curves are analytically equivalent. 
One easily verifies that any  analytically equivalent curves
share the same numerical semigroup.   {\bf For us, this  is the main point
of studying numerical semigroups.}    

By way of orienting ourselves to the subject, it may help to
state a few elementary results.  Below, let $S$ denote the
semigroup of a curve $c$.   

Fact.   $m(S) =1 \iff S = \N \iff c'(0) \ne 0 \iff c $ is equivalent to the line $t \mapsto (t,0, \ldots , 0)$.

 Fact:   $m(S) = 2 \iff  S = <2, 2k+1>$
 for some integer $k$  $\iff c$ is equivalent to the
 curve christened as the ``$A_k$ singularity'':   $(t^2, t^{2k+1})$, or, if $d> 2$,  $(t^2, t^{2k+1}, 0, \ldots, 0)$.
 
Fact: If $m(S) = m > 1$ and $n_1$ is the smallest generator of
$S$ besides $m$ then  $c$ is equivalent to a curve of the form  
$x_1 (t) = t^m,  x_2 (t) = t^{n_1} + \ldots $ where the `` $ \ldots$''
means term  of order greater than $n_1$ and 
where $ord(x_i (t)) > n_1$ for $i > 2$ , if $d> 2$.  

 \subsection{On well-parameterizing and real curves}

Consider  curves $c, \gamma: (\C, 0) \to (\C^d, 0)$
as being ``the same curve'' if there is a non-constant analytic map $\phi: (\C,0) \to (\C, 0)$
such that either  $c = \gamma \circ \phi$ or $\gamma = c \circ \phi$ holds. 
Curves that are ``the same''  in this sense  share the same image: 
$c(D) = \gamma(D')$ for $D, D' \subset \C$ sufficiently small appropriate  \nbhds of zero.

Given a curve $c$ which is not well-parameterized  we can always find another curve $\gamma$
that is ``the same''  as $c$ and  which is  well-parameterized.
To give a quintessential example, take $d =2$
and $c(t) = (x(t), y(t)) = (t^4, t^6)$.   Then $c$ is badly parameterized.
The curve 
$\gamma(t) = (t^2, t^3)$ is ``the same curve''  as $c$ since  $c = \gamma \circ \phi$
with  $\phi(t) = t^2$ and  $\gamma$ is well-parameterized.

Here is an alternative definition of ``well-parameterized''.  
A badly parameterized curve $c(t)$ can be factored as
$c = \gamma \circ \phi$ where $\phi:  (\C,0) \to (\C, 0)$ 
is a reparameterization having $\phi' (0) = 0$.  A ``well-parameterized curve' is a curve that   is not badly parameterized. 

Our original definition of   ``well-parameterized'' in terms of being one-to-one
does not work for curves over $\R$ because of maps like   $\phi(t) = t^3$
which are one-to-one over $\R$.  If $\gamma(t) = (t^2, t^3)$ 
then $\gamma \circ \phi$ is still one-to-one when viewed as a real curve
but we do   not consider it to be well-parameterized.

\subsection{Kunz Cone. }   The Kunz cone for  classifying multiplicity $m$
semigroups is a convex  polyhedral cone in $\R^{m-1}$ which provides 
a   direct and concrete way to parameterize all semigroups of multiplicity $m$.
The semigroups arise as  a subset of the lattice points in this cone. See \cite{Wilf} or    \cite{Oneill1}  and references therein.

To describe the  Kunz cone for multiplicity $m = 4$  we use 3 coordinates
$(x_1, x_2, x_3) \in \R^3$.   To obtain a  
semigroup associated to   a point  in this cone we    insist that the $x_i$ are positive  integers 
greater than $4$ and that 
$$x_i \equiv i  (mod 4).$$
\begin{definition}  An Apery  point in $\R^3$
is a point $(x_1, x_2, x_3) \in \R^3$ with $x_i \in \N$, $x_i > 4$
and $x_i \equiv i  $.
 Here, this and all subsequent   congruences ``$\equiv$'' are congruences mod $4$.
\end{definition}
The  semigroup associated to  an Apery   point is  $$S = <4, x_1, x_2, x_3>$$
It is essential here that we do not impose  the ordering   $x_1 < x_2 < x_3$.
It is also essential that the list $4, x_1, x_2, x_3$ need not be the list of minimal generators.
 
The Kunz cone is defined by the four  inequalities

\begin{equation}
    \text{Kunz cone} =
    \begin{cases}
     x_1 + x_2  &\ge x_3 \\
      x_3 + x_2 &\ge x_1 \\
      2x_1 &\ge x_2 \\
      2 x_3 & \ge x_2
    \end{cases}
    \label{Kunzcone}
  \end{equation}
  
To understand how these inequalities  arise,   look at the first one.
If $x_3 \ge  x_1 + x_2$ then, since $1+2 \equiv 3  (mod 4)$
we have that  $x_3 = x_1 + x_2 + 4k$
for some integer $k \ge 0$.  It follows that
$x_3 \in <4, x_1, x_2>$ or  that $<4, x_1, x_2> = <4, x_1, x_2, x_3>$.
So all   Apery points with $x_3 \ge x_1 + x_2$    yield the same semigroup,
and, notably, a semigroup  whose embedding dimension is $3$ (or perhaps even $2$) , not $4$.  Since
we want a bijection between points and semigroups, we exclude all Apery 
points having  $x_3 > x_1 + x_2$.  We keep the single point   $x_3 = x_1 + x_2$
so as to represent the semigroup $<4, x_1, x_2 >$. 

\begin{proposition}  The Apery points in the Kunz cone for multiplicity $4$
are in bijection with numerical semigroups of multiplicity $4$.
The points lying in the interior of the cone (i.e. all inequalities are strict) correspond to those
semigroups whose embedding dimension is $4$.  The Apery points lying
in the interior of the codimension $1$ faces of the cone correspond
to semigroups whose embedding dimension is  $3$ and those lying
on the rays of the cone (codimension $2$ faces) correspond to
semigroups of embedding dimension $2$.
\end{proposition}

Now, return to our problem of characterizing the multiplicity $4$ semigroups
whose minimal embedding dimension is $3$.  All the points on the
faces and rays of the Kunz cone have embedding dimension $3$ or less,
and so minimal embedding dimension $2$ or $3$.     In this way we
have  reduced our problem 
to the interior of the Kunz cone, which is to say, to those semigroups whose embedding dimension is
$4$.     
 
 \begin{theorem}  Let $S = <4, x_1, x_2, x_3>$
 be a  multiplicity 4 semigroup with embedding dimension 4,
 so that all the Kunz inequalities are strict.    Then,
 $S$ has minimal embedding dimension 3 if and only if
 \begin{itemize} 
 \item {a) } $max(x_1, x_2, x_3) \ne x_2$
 \item {b) } $max(x_1, x_2, x_3) > 2 min(x_1, x_2, x_3)$
  \end{itemize}
  \label{thm: main}
  \end{theorem}
  
  We have indicated these semigroups in figure XX.
  We are endebted here to the thesis of O'Sullivan \cite{OSullivan}  for 
  pointers in how to depict the Kunz cone via a slice.  Here we use the same
  slice $x_1 + x_2 + x_3 = 1$ that she used.  
  See figure \ref{kunzkite}.

  It is  easier  to understand the theorem and give its  proof if, instead of stating it as above, 
  we  run  through all $6$ possible orderings of the  $x_i$:
  
  (1)  $x_1 < x_2 < x_3$
  
  (2)  $x_2 < x_1 < x_3$
  
  (3) $x_1 < x_3 < x_2$
  
  (4)  $x_3 < x_1 < x_2$
  
  (5)  $x_2 < x_3 < x_1$
  
  (6)  $x_3 < x_2 < x_1$
  
  Then item (a) of the theorem says  that the orderings (3) and (4) do not  occur for 
  any semigroup  whose minimal embedding dimension  is 3 and  embedding dimension is $4$.   So, for example, the semigroup
  $<4, 7,9, 10>$ can only arise as the semigroup of a curve in 4 dimensional space.  
  
   Item (b) of the
  theorem asserts that if the $x_i$ satisfy  one of the orderings   (1), (2),  or (5), (6) 
  then the    inequality of (b) is  the necessary and sufficient conditions
  for  $<4, x_1, x_2, x_3>$ to have  minimal embedding dimension $3$
  and  embedding dimension $4$.  For example, if our generators are in the order of condition (1) and satisfy     $x_3 > 2x_1$
(and neccessarily $x_3 < x_1 + x_2$) then the  semigroup $<4, x_1, x_2, x_3>$ has  minimal embedding dimension $3$
and embedding dimension $4$.  An example 
of such a semigroup is  $<4, 5, 10, 11>$.
 
\section{proof}

  
 Let $n_1 < n_2 < n_3$ 
 be the list $x_1, x_2, x_3$ {\it permuted so as to be in numerical order}.
(For example, in case (5)   $n_1 =x_2, n_2 = x_3,  n_3 = x_1$.)
  Then, after a local diffeomorphism and  reparameterization  we can put our  curve
into the form
\beq
x = t^4,   y = t^{n_1} + a t^{n_1 + s_1} + \ldots,  z = t^{n_2} +  \ldots
\label{model form}
\eeq
Here  $s_1 > 0$ and the $\ldots$ in
both expansions means terms of order higher than the smallest written down,
so, in the case of $z$, of order higher than $n_2$.    
When written in this form we have $ord(x(t)) = 4,  ord(y(t)) = n_1,  ord(z(t)) = n_2$.
Now $ord(x^k (a y + bz)) \in <4, n_1, n_2>$
so that, in order to realize the largest generator $n_3$  of $S$
we require polynomials  which contain terms  quadratic in $y, z$.
Since  $ord( a y^2 + byz  + cz^2) \ge 2n_1$
we require
\beq
n_3 \ge 2n_1
\label{quad}
\eeq
if $S$ is to be the semigroup of a space curve.

We use inequality (\ref{quad}) to eliminate orderings (3) and (4).    In
both these cases $n_3 = x_2$. 
If we are  case (3) then  $n_1 = x_1$, 
but the strict Kunz inequality requires that   $x_2 < 2x_1$ or $n_3 < 2n_1$
violating  inequality (\ref{quad}).  Similarly in case (4) we have $n_1 = x_3$
and the strict Kunz inequality asserts again that $n_3 < 2n_1$, 
violating inequality (\ref{quad}).

To finish off the proof we exhibit space curves for the   other  four orderings  (1), (2) and (5) , (6)
  which realizes the given semigroup.  In   these cases
either $n_1$ or $n_2$ is $x_2$ so  that
$$n_1 + n_2 \equiv n_3$$
(Recall that  all   congruences ``$\equiv$'' are congruences mod $4$.) The strict Kunz inequalities then read
$$n_3 < n_1 + n_2 = ord(y(t) z(t)).$$
 We  are left with the fact that if we can find a polynomial
 whose order is $n_3$ then we must have that   $$max\{ n_2, 2 n_1 \}  <  n_3 < n_1 + n_2.$$
We are to show that   any  $n_3$ in this range
 and  congruent to $1$ or $3$ as appropriate,  can be realized 
 as the order of a polynomial pulled back to the curve.  
 We can set    
 $$n_3 = 2n_1 + s  \text{ for } 1 \le s < (n_2 -n_1) $$ and further fix the normal form  (\ref{model form})
 of the curve so that
 $$x = t^4, y = t^{n_1} + t^{n_1 + s}, z = t^{n_2}.$$
 Then  
 $$ y^2 =  t^{2n_1} + 2  t^{2 n_1 + s} + t^{2n_1 + 2s}$$
 If $n_1 =x_2$ then $2n_1 \equiv 0$ so that $2n_1 = 4k$.
 Then
 $y^2 - x^k = 2 t^{2n_1 + s} + \ldots$ so that
 $ord(y^2 -x^k) = n_3$.
 This takes care of cases (2) and (5).
 In the remaining two cases, (1) and (6),  $2n_1 \equiv 2$ and
 $n_2 \equiv 2$ so that 
 so that $2n_1 = n_2 + 4k$ which is
 the order of $x^k z$.  Then
 $y^2 -x^k z = 2 t^{2n_1 +s} + \ldots$
 and $ord(y^2 - x^k z) = n_3$.
 This takes care of all cases and completes the proof.
 
    QED
    
    \section{Appendix.   Teissier deciphered.}

We  describe   necessary and sufficient conditions
for a semigroup S to have  minimal embedding dimension  $me(S) = 2$.  
It  is a rewording  and   clarification of proposition 3.2.1    on p.  132-3  of  Teissier \cite{Teissier}.
 We  sketch a proof  
based on  the  Puiseux characteristic of a plane curve
and a recursion relation between    this  characteristic and the  semigroup of
the curve.  To do this, we will recall the Puiseux characteristic. 
 
 If $me(S) > 2$ then one implication  in  Teissier's    proposition still holds.
 This direction gives  a sufficient condition on  a semigroup $S$ 
 to have   $me(S) = d$ for   $d \le e$.   We describe what
we have been able to understand of  this part of the proposition.  

 \subsection{Divisor and factor vector of an increasing list}

Consider an increasing list of positive integers of length $e =g+1$:
$$\vec b = [b _0; b
_1,...,b _g];  b_i < b_{i+1}$$
Associated to this vector we have another integer
vector which we will call its   ``divisor vector'':
$$\vec e= [e_0, e_1, e_2, \ldots , e_g]$$
defined by the iteration scheme   
\beq 
e_i = {\rm g.c.d.}(b_0,..., b_i) \qquad, e_0 = b_0
\Leq{divisor list}
Note  that $e_i = {\rm g.c.d.}(e_{i+1}, b_i)$
and that all the $e_i$ are factors of $b_0$.   (Some of them may be $1$.)

The list $e_i$ is non-increasing and satisfies 
 $e_{i} | e_{i-1}$. It follows that for each $i > 0$    there is an integer $n_i  \ge  1$  given by  
$$n_{i}   = e_{i-1}/e_i  .$$
The list of $n$'s supply another    integer vector 
$$\vec n = [n_1, \ldots , n_g], $$
now  of length $g$, associated to $\vec b$. 
We call $\vec n$ the   ``factor vector''. Note that  $$n_g = e_{g-1}, $$
$$e_i = n_{i+1} n_{i+2} \ldots n_g.$$
In particular
$$b_0 = n_1 n_2 \ldots n_g.$$

\subsection{Planar semigroups lie on the boundary of the Kunz cone} 
Set $$m = b_0$$
in order to  denote  the multiplicity of the curve corresponding
to the semigroup $S$.
Since each $n_i \ge 2$ we have that $m \ge 2^g$
or $log_2 m \ge g$.  Now, as long as $m \ge 6$
we have that $m > log_2 m  + 2$
so  $m-1 > g+1$.  Whem $m = 3, 5$ we have $g=1$ or  $e=2$ and so the planar semigroups
lie on the 1-dimensional faces of the Kunz cone.  
In the case   $m = 4$ we just dealt with we can  have $g=2$ and
so $e =3 = m-1$ and the planar semigroups lie on the faces of the Kunz cone.  
It follows that in all situations in which $m > 2$
the planar semigroups lie on the boundary of the
Kunz cone.

\subsection{Conditions for a semigroup to be that of a plane curve}

\begin{proposition}  
Let $S = <b_0, b_1, \ldots , b_g>$  be a numerical semigroup given by
its minimal generators $b_i$ listed in order.  
Write  $\vec e = [e_0, e_1, \ldots , e_g]$ for   the 
divisor vector of $[b_0, b_1, \ldots , b_g]$ and   $\vec n = [n_1, n_2, \ldots , n_g]$  
 for its  factor vector  ($n_i = e_{i-1}/e_i$).   
 
 $S$ has   minimal embedding dimension $2$  if and only if  the following
 three conditions hold for its associated three  integer vectors
 
 \noindent (1):   the divisor vector  $\vec e$   is strictly decreasing and ends with $1$.  
 
\noindent (2):  for      $i =1, \ldots, g$ we have 
 $n_i b_i \in < b_0, b_1,  \ldots , b_{i-1} >$

\noindent (3):  for      $i =1, \ldots, g-1$ we have $n_i b_i < b_{i+1}$
\label{main2}
\end{proposition}   

We sketch a proof of   the proposition at the
very end of this appendix.  We do so  by translating 
the proposition, and the planar semigroups, into the languague of  Puiseux characteristics.
We recall that language in time for the proof at the end.

\begin{example}  If  $S = <b_0, b_1, \ldots >$ and  $b_0, b_1$ are relatively prime then $e_1 = 1$.
 If condition (1) holds then there cannot be an $e_2$,
 so  it must be that
 $g=1$.  Thus   $S = <b_0, b_1>$ if $S$ comes from a planar semigroup whose
 first two generators are relatively prime.   These are represented by the curves 
 discussed in the first paragraph of this paper. 
\end{example}

\begin{corollary}
If $m = b_0$ is the multiplicity of a semigroup $S = <b_0, b_1, \ldots , b_g >$
given by its  minimal generating set $b_i$ and if $S$ is planar, i.e. has  minimal embedding dimension $2$,
 then the embedding dimension $g+1$ of
$S$ is less than or equal to $1$ plus the number of prime factors of
$m$ counted with multiplicity.
\end{corollary} 

{\sc Proof of Corollary.} The longest we can make the divisor vector $[e_0, e_1, \ldots , e_g]$
and keep  it strictly decreasing  
beginning with $m = e_0$  is $1$ plus the number of prime factors of $m$.
We do this by omitting one factor at a time  from $m$so that $e_i = e_{i-1}/p_i$
where the $p_i$ exhaust the  prime factors, taken with multiplicity,  of $m$. QED

\begin{remark}   Condition (1) is extraneous in that it is implied
by condition (2) and the   assumption that the $b_i$ form a minimal set of generators.
For, if (2) holds, then the non-increasing list $e_i$ has to be strictly decreasing.
To see this we  argue by contraposition. If $e_i$ is not strictly decreasing then
there is an $i$ such that $e_{i-1} = e_i$ in which case $n_i =1$.
But then (2) says that $b_i \in < b_1, \ldots, b_{i-1}>$
which implies that the generator $b_i$ can be omitted from our  list of generators
and we still have a generating set for $S$.  
\end{remark}

\subsection{An aside: self-duality and planarity}

The conductor $c$ of a numerical semigroup $S$ is the smallest element of  $S$
such that all integers greater than $c$ lie in $S$.    The set of gaps
of a numerical semigroup is the finite set $\N \setminus S$.
A numerical semigroup is called ``self-dual'' if  the cardinality of the set of its gaps 
is half of its  conductor $c$.   
It is a well-known fact that planar semigroups are self-dual.
 Here is a  notable example, described as
``3.2.3 Remark'' on p. 133 of Teissier, of a self-dual semigroup which is not
that of a planar curve.

\begin{example}  [Teissier]
Set $S = <9, 21, 22>$.   A session with GAP yields that
the conductor of $S$ is 78 and the number of gaps is 39
so that $S$ is self-dual.
The divisor list of $S$ is $\vec e = [9, 3, 1]$
and the factor list is $\vec n = [3, 3]$.   
We have 
$(b_0, b_1, b_2) = (9,21, 22)$.  Since
$n_1 b_1 = 63 > b_2 = 22$ this semigroup fails condition
of  (3) 
of the proposition so cannot be the semigroup of a planar curve.
\end{example}  

\subsection{Teissier's proposition for   non-planar curves}

Teissier wrote  his proposition so that one direction  of its implications
holds  for  any minimal embedding dimension $d$ in place of $d=2$.  
I describe the statement as   I understand it.  Take it with a bit of scepticism,
since I do not understand his proof.

Suppose that a   semigroup $S$ is  given by its  minimal generators 
$b_i$ as before.  Drop condition (1) of the proposition.
Suppose that condition (2) holds.  (Teissier labels this condition (1)).
And suppose that condition (3) holds,  but not for all $g-1$ listed indices, 
but instead for  
$\ell = g+1  - d$ of these indices.
Then there exists  an analytic curve $\gamma:  (\C, 0) \to (\C^d, 0)$
having $S$ as its semigroup.  

He  says nothing about the converse direction of the implication.
It is generally false as examples with $e = me$ and $b_0, b_1$
relatively prime show.

 Within the proof on page 133 Teissier
says that condition (2) implies that the canonical curve
$x_i = t^{b_i}$ is a complete intersection.  Is it a necessary and
sufficient condition to be a complete intersection? 
I do not know.

\subsection{Strategy of proof of the Proposition}

We show how to go from a plane curve to its Puiseux characteristic. 
We describe a transformation taking us  from the Puiseux characteristic to the
semigroup of the curve.  It is   straightforward to check that
the resulting semigroup satisfies the proposition.  Conversely,
starting from the semigroup we can   reverse the Puiseux-to-semigroup transformation to recover the Puiseux
characteristic  from the semigroup, provided the semigroup  satisfies the conditions  of
the proposition.  And given the  Puiseux characteristic it is immediate
to write down a plane curve with this Puiseux characteristic. 
In order to implement this strategy we begin with reviewing the Puiseux characteristic of a plane curve. 

\subsection{The Puiseux characteristic  of a plane curve}

The \Pc \ of an analytic plane curve germ $c: (\mathbb C,0)\to
\mathbb C^2$ is an increasing  vector $\vec \lambda = [\lambda_0, \lambda_1, \ldots,  \lambda_g]$
of positive integers encoding the key exponents which
arise in the power series expanison of $c(t) = (x(t), y(t))$.
  A reparameterization and change of coordinates
  puts the curve into the form  
\begin{equation}
\label{eq-P1} x = t^{m}, \ \ y =  a_n t^{n} +  a_{n+1}t^{n+1} + a_{n+2}  t^{n+2} + \cdots , \ \  
\end{equation}
where $1 < m < n$ and $a_n \ne 0$.
Then
$$\lambda_0 = m.$$
 
We can assume that $n$ is not a multiple of $m$ in the expansion of the curve.  For if
$n = km$ then the diffeomorphism $(x, y) \mapsto (x, y- a_n x^k)$
kills the term $a_n t^n$.    We can kill all powers of $t^m$
arising in the power series of $y(t)$ by this same trick.
With this in mind,  write $supp(y) = A$ for the set of all exponents occuring
in the power series of $y$.  Thus
$$A = \{ i \in \N:  a_i \ne 0 \} \text{ where }  y(t) = \Sigma a_i t^i$$
We have just seen that we can get rid of all  multiple of $m$
arising in $A$ by applying a   diffeomorphism $(x, y) \mapsto (x, y- f(x))$.
 We have ``seived out'' $m$ from $A$.  
\begin{definition} If $A \subset \N$ and $m \in \N$ then the
m-seive of $A$, denoted $[m; A]$  is the set $A$ with all multiples of 
$m$ deleted.  In set notation, $[m;A] =  A \setminus m \N$.  \end{definition}
Since $c(t)$ is well-parameterized, the m-seive  of   $A$ is not empty.

Let $\lambda_1$ be the smallest element of $[m; A]$.  
Thus, $\lambda_1 = n$ above, assuming that $m$ does not divide $n$.
 If $gcd(m, n_1) =1$
the process ends and the Puiseux expansion of $c(t)$ is $[m, \lambda_1]$.
Otherwise, write $e_1 = gcd(m, \lambda_1)$.

We continue by setting $A_2 = [e_1; A]$ and choosing $\lambda_2$ to be the smallest element
of  $A_2$.    We write $e_2 = gcd(e_1, \lambda_2)$.  If $e_2 = 1$ we stop
and Puiseux characteristic is $[m, \lambda_1, \lambda_2]$.  Otherwise we set $A_3 = [e_2; A_2]$
and take $\lambda_3$ to be the smallest element of $A_3$.  
Iterate  the process:   $\lambda_i = min(A_i)$,  $e_i = gcd(e_{i-1}, \lambda_i)$, and 
$A_{i+1} = [e_i; A_i]$.  We stop with $\lambda_g$ when    $e_g = 1$.
We have constructed  an increasing list
$m = \lambda_0  < \lambda_1 < \lambda_2 < \ldots < \lambda_g$
of integers which have no common divisor.    Our Puiseux expansion is this list.

In the process, we have also constructed the divisor vector
of $\vec n$, namely the $e_i$.  It is a      strictly decreasing vector $e_i < e_{i-1}$.

\begin{example} The Puiseux characteristic of  the plane curve germ
$$(t^8, \ t^{16} + t^{20} + a_{22}t^{22} + a_{26}t^{26} +
a_{27}t^{27})$$ is $[8; 20, 22, 27]$ provided that $a_{22},
a_{27}\ne 0$.   The divisor vector of this Puiseux characteristic
is $\vec e = [8, 4, 2, 1]$.   Its factor vector is $[2,2,2]$.
\label{eg 8...} 
\end{example}

\subsection{ Characterizing Puiseux characteristics }  Here are the necessary and sufficient conditions for 
a list of integers, labelled now $[\lambda_0, \lambda_1, \ldots, \lambda_g]$ to be the Puiseux characteristic
of some curve.  

\medskip

\noindent 1.  $\lambda_0 > 1$ and the list is  strictly increasing: $\lambda_i < \lambda_{i+1}$. 
\medskip

\noindent 2.  Its associated divisor vector   $\vec e = [\lambda_0,  e_1, \ldots , e_g]$  is strictly decreasing and
ends with $1$.

\medskip

\begin{rmk} A  \Pc of length 2 consists of a pair  of relatively prime integers
$[\lambda_0, \lambda_1]$ with $\lambda_0 < \lambda_1$.  
By fiat (or a logical contortion if you prefer) the only \Pc of length 
$1$ is  the vector $\lambda = [1]$.  
\label{rmk: rel prime}
\end{rmk}
  
\subsection{Puiseux to Semigroup} 
The following  iteration scheme  recovers the semigroup generators
$<b_0, b_1, \ldots, b_g>$  of the semigroup of a plane curve  from
its  Puiseux characteristic $[\lambda_0; \lambda_1, \ldots , \lambda_g]$.
I am endebted to Teissier, 2.2.1, formula (*) p. 122 for this scheme.  
   C.T.C. Wall 's book , Prop. 4.3.8 on p. 86 and  five or so pages on
   either side also covers this recursion formula.  
{\it We will see how this process works for
the curve of example \ref{eg 8...} at the   end of this subsection.}  

\medskip 
$$b_0 = \lambda_0$$

$$b_1 = \lambda_1$$

$$b_2 = \lambda_2 -\lambda_1 + n_1 \lambda_1$$

And inductively  
$$b_i = \lambda_i - \lambda_{i-1} + n_{i-1} b_{i-1}$$

Using $n_{i-1} > 1$ one easily verifies by induction that
$$b_i > \lambda_i,  i >1.$$
and that $b_i \in < \lambda_1, \lambda_2, \ldots,  \lambda_i >$.

\medskip

I will not rederive the recursion  formula but simply content myself with the understanding the generator $b_2$.
So suppose that $e_1 = gcd (m , \lambda_1) > 1$ where $\lambda_0 = m$.
Then we can put $c(t)$ into the form
$$x = t^m,  y = t^{\lambda_1} + a t^{\lambda_2} + \ldots;  a \ne 0.$$
We have that $n_1 e_1 = m$ and $\lambda_1 = \beta e_1$ for some integer $\beta$.
Then
$$y^{n_1} = t^{n_1 \lambda_1} + n_1 a t^{\lambda_2} t^{(n_1 -1) \lambda_1}  +  O(t^j)  , j > \lambda_2 + (n_1 -1) \lambda_1  $$
But $n_1 \lambda_1 = \beta n_1 e_1 = \beta m$ so that
$$x^{\beta}  = t^{n_1 \lambda_1}$$
which shows that $ord(x^{\beta} - y^{n_1}) = \lambda_2 + (n_1 -1) \lambda_1$.
Hence $\lambda_2 + (n_1 -1) \lambda_1 \in S$, where $S$ is the semigroup of this curve. 
It is not difficult to show that this integer is the smallest   element of 
$S$ not lying in $<m, \lambda_1> \subset S$ and
hence this integer is the  next generator $b_2$ of $S$ after $b_1 = \lambda_1$.

\begin{rmk}The $b_i$ obtained from this formula have  the same divisor vector
$\vec e = [e_0, e_1, \ldots , e_g]$ and factor vector $[n_1, n_2, \ldots , n_g]$
as that  of the Puiseux vector.  This follows directly from the recursion relation.  
  Using the same recursion relation, we can
invert an associated matrix and solve for the $\lambda_i$ given the $b_i$.  
\label{rmk: same div list}
\end{rmk}

\begin{example}   Return to our earlier curve
$(t^8, \ t^{16} + t^{20} + a_{22}t^{22} + a_{26}t^{26} +
a_{27}t^{27})$ whose Puiseux characteristic is $[8; 20, 22, 27]$ provided that $a_{22},
a_{27}\ne 0$.   To find its semigroup compute that
$$[e_0, e_1, e_2, e_3] = [8, 4, 2, 1]$$
while
$$[n_1, n_2, n_3] = [2, 2,2]$$
We compute $b_0 = \lambda_0 = 8,  b_1 = \lambda_1 = 20$
$$b_2 = 22 - 20 + 2* 20 = 42$$
$$b_3 = (27-22) + 2*b_2 = 5+ 84 = 89$$
so that $S = <8, 20, 42, 89 >$

As  a reality check   we verify that   we can get the integers
$42, 89$ as orders of polynomials pulled back to the curve.  Take the 
case $a_{22} = a_{27} = 1,  a_{26} = 0$ for simplicity so that
$y(t) = t^{20} + t^{22} + t^{27}$.  Then  
$y^2 - x^5 = (t^{40} + 2 t^{42} + 2 t^{47} + \ldots .) - t^{40} = 2 t^{42} + 2 t^{47}$
which has order $42$.   To get $89$ note that  
$(y^2 - x^5)^2 = 4 t^{84} + 8 t^{89} + ...$
Now $84$ is a multiple of $4$.   Since $[8, 20] = 4[2, 5]$
and since the conductor of $<2, 5>$ is $4$, from $16$ onwards all integers which
are  
multiples of $4$ occur  in any  semigroup having $8$ and $20$ as generators.   
It follows that we can kill the $t^{84}$ term of $(y^2 - x^5)^2$ with some polynomial in $x$, and $y$.
Indeed $84 = 64 +20$ is the valuation of  $x^8 y$.  We see that 
$$(y^2 - x^5)^2 - 4 x^8 y  = 8 t^{89} + ...$$
which has a valuation of $89$. 

\end{example}

\subsection{Sketch of a proof of Teissier's proposition \ref{main2}}

Given  a plane curve, form its Puiseux expansion.
Run the recursion to obtain the minimum generating
list for the   semigroup of the curve.
Use remark \ref{rmk: same div list} to see that the divisor vector
of the list semigroup generators is strictly decreasing.  Inductively verify
conditions (2) and (3).

Conversely, given a semigroup whose minimal generating list satisfies
the given relation run the recursion relation backwards to obtain a
  Puiseux expansion $[\lambda_0, \lambda_1, \ldots, \lambda_g]$ whose
  semigroup is the give semigroup.  
The curve $x = t^{\lambda_0},  y = t^{\lambda_1} +  t^{\lambda_2} + \ldots +  t^{\lambda_g}$
has $S$ as its semigroup.

\section*{Acknowledgement}
I'd like to thank Justin Lake a graduate student finishing up at  UCSC, the algebraic geometers Gary Kennedy and Lee McKewan 
of Ohio, and the numerical semigroup specialist Chris O'Neill of San Diego
for conversations during the research and writing of this piece.  I would like to thank Emily O'Sullivan 
for redrawing my rough sketch of figure 1, the Kunz Kite.
 
\end{document}